\documentclass[11pt]{article}
\usepackage{amsmath}[1996/11/01]
\usepackage{amssymb,amsthm,amsfonts,latexsym,psfig,epsfig,here}

\usepackage{graphics,psfig,epsfig}
\usepackage[colorlinks=true,
linkcolor=webgreen,
filecolor=webbrown,
citecolor=webgreen]{hyperref}
\definecolor{webgreen}{rgb}{0,.5,0}
\definecolor{webbrown}{rgb}{.6,0,0}

\newtheorem{theorem}{Theorem}[section] 

\newtheorem{corollary}[theorem]{Corollary}
\newtheorem{conjecture}[theorem]{Conjecture}

\newtheorem{definition}[theorem]{Definition}
\newtheorem{remark}[theorem]{Remark}

\newcommand{\FF}{\mathbb F}

\newcommand{\eqn}[1]{(\ref{#1})}

\begin{document}

\begin{center}
{\Large {\bf On Single-Deletion-Correcting Codes}} \\
\vspace{1.5\baselineskip}
{\em N. J. A. Sloane} \\
\vspace*{1\baselineskip}
Information Sciences Research, AT\&T Shannon Labs \\
180 Park Avenue, Florham Park, NJ 07932-0971, U.S.A. \\
\vspace{1.5\baselineskip}
{\small\rm Dedicated to Dijen Ray-Chaudhuri on the occasion of his 65th birthday}
\vspace{1.5\baselineskip}
\end{center}


\begin{abstract}
This paper gives a brief survey of binary single-deletion-correcting codes. The 
Varshamov-Tenengolts codes {\em appear} to be optimal, but many 
interesting unsolved problems remain. 
The connections with shift-register sequences also 
remain somewhat mysterious. 
\end{abstract}


\section{Introduction}\label{sectI}
The possibility of packet loss on internet transmissions has renewed interest in deletion-correcting codes.
(Of course there are many other applications of such
codes, including magnetic recording, although in that case there are usually additional conditions that must be satisfied.)
This paper considers the very simplest family of such codes, binary block codes capable of correcting single deletions.
Even for these codes there remain several apparently unsolved problems.

It is surprising, but these codes do not appear to be surveyed in any of the usual references
(\cite{MS77}, \cite{PH98}, etc.).
This paper is a first attempt at such a survey.

Proofs are given of a number of results, either because the new proofs are simpler or because the original sources are hard to locate\footnote{and when located are sometimes poorly translated or badly photocopied!}.

\begin{definition}\label{de1}
For a vector $u \in \FF_q^n$, let $D_e (u)$ denote the set of $e$-th order descendants, i.e. the set of vectors $v \in \FF_q^{n-e}$ that are obtained if $e$ components are deleted from $u$.
A subset $C \subseteq \FF_q^n$ is said to be an {\em $e$-deletion-correcting code} if $D_e ( u) \cap D_e (v) = \emptyset$ for all $u,v \in C$, $u \neq v$.
Our problem is to find the largest such code.
In this paper we mostly consider the simplest case, $q=2$ and $e=1$.
\end{definition}

The {\em deletion distance} $dd(u,v)$ between vectors $u,v \in \FF_q^n$ is defined to be one-half of the smallest number of deletions and insertions needed to change $u$ to $v$.
Then $C$ is $e$-deletion-correcting if and only if $dd(u,v) \ge e+1$ for $u,v \in C$, $u \neq v$.
(For $dd(u,v) \le e$ if and only if there is a vector $x$ that can be reached from $u$ by at most $e$ deletions and also from $v$ by at most $e$ deletions,
and then $C$ cannot correct $e$ deletions.)

Consider the graph $G_n$ having a node for every vector $u \in \FF_q^n$, with an edge joining the nodes corresponding to $u,v \in \FF_q^n$, $u \neq v$, if
and only if $v$ can be obtained from $u$ by a single deletion and insertion, i.e. if and only if $D_1 (u) \cap D_1 (v) \neq \emptyset$.
The deletion distance $dd(u,v)$ is the length of the shortest path from $u$ to $v$
(this shows that $dd$ is indeed a metric).

In particular, a single-deletion-correcting code corresponds to an independent set in $G_n$.
One can now attempt to calculate the sizes of the largest independent sets by computer.
In the binary case we find that the largest single-deletion-correction codes of lengths $1,2, \ldots, 8$ have sizes
\begin{equation}\label{EqI1}
1,2,2,4,6,10,16, \ge 30 ~.
\end{equation}
The last entry in (\ref{EqI1}) was kindly computed by my colleague David Johnson.
Unfortunately $G_8$ is too large for present computers and 30 is at present only a lower bound on the size of a maximal independent set.\footnote{Postscript: David Applegate has since used CPLEX's integer programming subroutines
(which combine ordinary linear programming with branch-and-bound)
to confirm that the largest single-deletion-correcting code of length 8 does indeed have size 30.}

However, \eqn{EqI1} turns out to be a useful hint.
When one looks up this sequence in \cite{EIS}, one finds a unique
matching sequence, number
\htmladdnormallink{A16}{http://www.research.att.com/cgi-bin/access.cgi/as/njas/sequences/eisA.cgi?Anum=000016},
whose initial terms $N_1$, $N_2$, $N_3, \ldots$ are
\begin{equation}\label{EqI2}
1,1,2,2,4,6,10,16,30,52,94,172, 316, 586, \ldots
\end{equation}
and whose $n$th term is given by
\begin{equation}\label{EqI3}
N_n = \frac{1}{2n} \sum_{{\rm odd} \, d|n} \phi (d) 2^{n/d}, ~n \ge 1 \,,
\end{equation}
where the sum is over all odd divisors $d$ of $n$ and $\phi$ is the Euler
totient function (sequence
\htmladdnormallink{A10}{http://www.research.att.com/cgi-bin/access.cgi/as/njas/sequences/eisA.cgi?Anum=000010}).
The references cited for sequence A16 indicate that it has arisen in connection with the
enumeration of shift-register sequences \cite{Go67} and
tournaments \cite{Br80}.
However there was (at that time) no reference to indicate that this sequence has any
connection with codes, nor was there any apparent connection between the 
shift-register sequences and deletion-correction codes.

More conventional search methods, in particular, consulting some well-known papers of Levenshtein \cite{Lev65}, \cite{Lev65a} on codes for correcting deletions, turned up many other relevant references.
Some of these will be discussed further in Section \ref{sectH}.
The most interesting codes are those of Varshamov and Tenengolts \cite{VT65}.
In \cite{VT65} they present a family of codes depending on a certain parameter $a$.
When $a$ is taken to be 0, these codes have size $N_{n-1}$ (see \eqn{EqI3}) and
thus match \eqn{EqI1}.
These codes are the subject of Section \ref{sectVT}.

Sections \ref{sectSR} and \ref{sectT} will discuss the connection with shift-registers and tournaments, and Section \ref{sectM} contains some
general remarks about the number of descendants of a vector.
The final section, Section \ref{sectH}, gives a brief discussion of other papers on deletion-correcting and related codes.

\section{The Varshamov-Tenengolts codes}\label{sectVT}
\begin{definition}\label{de21}
For $0 \le a \le n$, the
{\em Varshamov-Tenengolts code} $VT_a (n)$ consists of all binary vectors
$(x_1, \ldots, x_n )$ satisfying
\begin{equation}\label{EqVT1}
\sum_{i=1}^n ix_i \equiv a ~(\bmod~n+1) \,,
\end{equation}
where the sum is evaluated as an ordinary rational integer.
\end{definition}

As will appear, the codes with $a=0$ contain the most codewords.
The first few such codes are
\begin{eqnarray}\label{EqVT1a}
VT_0 (1) & = & \{0\} \nonumber \\
VT_0 (2) & = & \{00, 11\} \nonumber \\
VT_0 (3) & = & \{000, 101\} \nonumber \\
VT_0 (4) & = & \{0000, 1001, 0110, 1111\} \nonumber \\
VT_0 (5) & = & \{ 00000, 10001, 01010, 110011, 11100, 00111 \} \,,
\end{eqnarray}
of sizes 1,2,2,4,6, matching \eqn{EqI1} and \eqn{EqI2}.
These codes were introduced in \cite{VT65} for correcting errors on a $Z$-channel (or asymmetric channel).
Similar constructions have been used in \cite{BR82} and also in
\cite{GS80} and \cite{Kl81} to construct constant weight codes.

Levenshtein \cite{Lev65}, \cite{Lev65a} observed that the Varshamov-Tenengolts codes could be used for correcting single deletions, proving this by giving the following elegant decoding algorithm.

\subsection*{Decoding algorithm}
\begin{itemize}
\item
Suppose a codeword $x = (x_1, \ldots, x_n ) \in VT_a (n)$ is transmitted, the symbol $s$ in position $p$ is deleted, and $x' = (x'_1, \ldots, x'_{n-1} )$ is received.
Let there be $L_0$ 0's and $L_1$ 1's to the left of $s$, and $R_0$ 0's and $R_1$
1's to the right of $s$ (with $p = 1+ L_0 + L_1 )$.
\item
We compute the weight $w = L_1 + R_1$ of $x'$ and the new checksum
$\sum_{i=1}^{n-1} ix'_i$.
If $s = 0$ the new checksum is $R_1$ $(\le w )$ less than it was before, and if $s=1$ it is $p+ R_1 = 1+ L_0 + L_1 + R_1 = 1 + w + L_0$
$(> w )$ less than it was before.
(These numbers are less than $n+1$ so there is no ambiguity.)
\item
So if the deficiency in the checksum is less than or equal to $w$ we know that a 0 was deleted, and we restore it just to the left of the rightmost $R_1$ 1's.
Otherwise a 1 was deleted
and we restore it just to the right of the leftmost $L_0$ 0's.
\end{itemize}
\begin{table}[htb]
\caption{Number of codewords in Varshamov-Tenengolts code $VT_a (n)$.}
$$
\begin{array}{c|ccccccccc}
n \setminus a & 0 & 1 & 2 & 3 & 4 & 5 & 6 & 7 & 8 \\ \hline
1 & 1 & 1 \\
2 & 2 & 1 & 1 \\
3 & 2 & 2 & 2 & 2 \\
4 & 4 & 3 & 3 & 3 & 3 \\
5 & 6 & 5 & 5 & 6 & 5 & 5 \\
6 & 10 & 9 & 9 & 9 & 9 & 9 & 9 \\
7 & 16 & 16 & 16 & 16 & 16 & 16 & 16 & 16 \\
8 & 30 & 28 & 28 & 29 & 28 & 28 & 29 & 28 & 28
\end{array}
$$
\label{TVT1}
\end{table}

The sizes
$|VT_a (n)|$ of the first few codes are shown in Table \ref{TVT1}.
(This array forms sequence
\htmladdnormallink{A53633}{http://www.research.att.com/cgi-bin/access.cgi/as/njas/sequences/eisA.cgi?Anum=053633}
in \cite{EIS}.)
These numbers were studied by Varshamov \cite{Var65} and Ginzburg \cite{Gi67}, but the following simple formula appears to be new.

\begin{theorem}\label{TH1}
\begin{equation}\label{EqVT2}
| VT_a (n) | = \frac{1}{2(n+1)} \sum_{d| n+1 \atop d~{\rm odd}}
\phi (d) \frac{\mu \left( \frac{d}{(d,a)} \right)}{\phi \left( \frac{d}{(d,a)} \right)} 2^{(n+1)/d} ~,
\end{equation}
where $\mu (n)$ is the M\"{o}bius function
$($\htmladdnormallink{A8683}{http://www.research.att.com/cgi-bin/access.cgi/as/njas/sequences/eisA.cgi?Anum=008683}$)$, and
$(d,a) = gcd (d,a)$.
\end{theorem}

\begin{proof}
Write $w_a (n) = | VT_a (n) |$.
We will calculate $w_a (n-1)$, assuming throughout that $n \ge 1$.
It follows from the definition of these codes that the generating function
$$f(z) = \sum_{a=0}^{n-1} w_a (n-1) z^a$$
is equal to
$$\prod_{k=1}^{n-1} (1+z^k ) ~\bmod~ z^n -1 ~.
$$
Let $\xi = e^{2\pi i /n}$.
Then
$$f(\xi^j ) = \sum_{a=0}^{n-1} w_a (n-1) \xi^{ja} =
\prod_{k=1}^{n-1} (1+ \xi^{jk} ) ,~
j=0,\ldots, n-1 \,.
$$
We solve this by taking an inverse discrete Fourier transform
(cf. \cite{Ko88}, Chap. 97) to obtain
$$w_a (n-1) = \frac{1}{n}
\sum_{j=0}^{n-1} f(\xi^j) \xi^{-ja} \,.
$$

Since
$$\prod_{k=0}^{n-1} (z - \xi^k ) = z^n -1 ~,$$
we can calculate $f(\xi^j )$ explicitly.
An elementary calculation gives
$$f(\xi^j ) = \left\{ \begin{array}{lll}
2^{g-1} & \mbox{if} & d=n/g ~\mbox{is odd} , \\ [+.1in]
0 & \mbox{if} & d=n/g ~\mbox{is even},
\end{array}
\right.
$$
where $g= gcd (n,j)$.
Therefore
$$w_a (n-1) = \frac{1}{2n} \sum_{d|n \atop d~{\rm odd}}
2^{n/d} \sum_{j=1 \atop gcd (n,j) = n/d}^n \xi^{-ja}
$$
which becomes, writing $j = kn/d$,
$$ = \frac{1}{2n}
\sum_{d|n \atop d~{\rm odd}} 2^{n/d} \sum_{k=1 \atop (k,d) =1}^d
e^{-2 \pi i ka /d} \,.
$$
The innermost sum is a Ramanujan sum $c_d (a)$
(\cite{Ap76}, p. 160), which simplifies to
$$c_d (a) = \phi (d) \frac{\mu \left( \frac{d}{(d,a)} \right)}{\phi \left( \frac{d}{(d,a)} \right)}
$$
(\cite{Ap76}, p. 164).
\end{proof}

\begin{corollary}\label{TH2}
\begin{eqnarray}\label{EqVT3}
\mbox{(i)}~ |VT_0 (n) | & = & \frac{1}{2(n+1)}
\sum_{d | n+1 \atop d~{\rm odd}} \phi (d) 2^{(n+1)/d} ~, \\
\label{EqVT4}
\mbox{(ii)} ~ |VT_1 (n) | & = & \frac{1}{2(n+1)}
\sum_{d| n+1 \atop d~{\rm odd}} \mu (d) 2^{(n+1)/d} ~,
\end{eqnarray}
\begin{itemize}
\item[(iii)]
For any $a$,
\begin{equation}\label{EqVT5}
|VT_0 (n) | \ge | VT_a (n) | \ge |VT_1 (a) | \,.
\end{equation}
\end{itemize}
\end{corollary}

\begin{remark}
{\rm (i) and the left-hand inequality in (iii) are due to Varshamov \cite{Var65},
and (ii) and the right-hand inequality in (iii) to Ginzburg \cite{Gi67}.}
\end{remark}

\begin{proof}
(i) and (ii) follow immediately from Theorem \ref{TH1},
as does the left-hand side of (iii) using $\mu (k) \le \phi (k)$ for all $k$.
To establish the right-hand side of (iii),
let $p$ be the smallest odd prime dividing both $n+1$ and $a$ (if no such prime exists then $|VT_a (n) | = |VT_1 (n)|$).
The terms in the expressions for $|VT_a (n) |$ and $|VT_1 (n) |$ agree for
$d < p$, and at $d=p$ the term in $|VT_a (n)|$ exceeds that in $|VT_1 (n) |$ by $p2^{n/p}$.
It is easy to check that the remaining terms can never make the sum in
$|VT_1 (n)|$ catch up with the sum in $|VT_a (n) |$.
\end{proof}

\section*{Optimality}
It is more difficult to obtain upper bounds for deletion-correcting codes than for conventional error-correcting codes, since the disjoint balls $D_e (u)$ associated with the codewords (see Section \ref{sectI}) do not all have the same
size.
Furthermore the metric space $(\FF_2^n, dd )$ is not an association scheme and so there is no obvious linear programming bound.

The size of $D_1 (u)$ is easily seen to be equal to $r(u)$, the number
of runs in $u$.
Furthermore the number of vectors in $\FF_2^n$ with $r$ runs is
$2 {\binom{n-1}{r-1}}$.
(We will discuss $|D_e (u) |$ further in Section \ref{sectM}.)

Let $A(n,e)$ denote the size of the largest $e$-deletion-correcting binary code of length $n$, and call a code $C$ {\em optimal} if $|C| = A(n,e)$.
The values of $A(n,1)$ for $n \le 9$
were given in Section \ref{sectI}, and show that $VT_0 (n)$ is optimal for $n \le 9$.

For large $n$, the codes $VT_0 (n)$ are certainly close to being optimal, since on the one hand we have
\begin{equation}\label{EqVT6}
|VT_0 (n) | \ge \frac{2^n}{n+1} ~,
\end{equation}
from \eqn{EqVT5}, and on the other hand we have the following result of
Levenshtein \cite{Lev65}:
\begin{theorem}[\cite{Lev65}]\label{TH3}
$$A(n,1) \sim \frac{2^n}{n}, \quad\mbox{as}\quad n \to \infty \,.
$$
\end{theorem}

\begin{proof}
\eqn{EqVT6} gives a lower bound.
Let $C$ be an optimal code.
Following Levenshtein, let $C_0$ denote the subset of $C$ consisting of the vectors
$u \in C$ with
$$\frac{n}{2} - \sqrt{n \log n} \le r(u) \le \frac{n}{2} +
\sqrt{n \log n}
$$
and let $C_1 = C \setminus C_0$.
Since the sets $D_1 (u)$, $u \in C$, must be disjoint,
$$|C_0| \le \frac{2^{n-1}}{\frac{n}{2} - \sqrt{n \log n}} \lesssim
\frac{2^n}{n} ~.
$$
Furthermore,
$$|C_1 | \le 2 \sum_{r=1}^{\frac{n}{2} - \sqrt{n \log n}} 2 {\binom{n-1}{r-1}} ~,
$$
which is much smaller than $2^n /n$.
\end{proof}

In a later paper, Levenshtein \cite{Lev92} defines a code $C$ to be {\em perfect} if the balls $D_e (u)$, $u \in C$, partition the set $\FF_2^{n-e}$.
In \cite{Lev92} he proves the remarkable
fact that {\em all} the codes $VT_0 (n)$, $VT_1 (n)$, $VT_2 (n) , \ldots$ are perfect single-deletion-correcting codes.
The argument, not reproduced here, is essentially just a refinement of the
decoding algorithm for these codes given above.

It is initially
surprising that perfect codes of the same length can have different numbers of codewords, but this is explained by the fact that the balls $D_1 (u)$ have different sizes.

In view of this and the result in \eqn{EqVT5}, it is tempting to make the following conjecture.

\begin{conjecture}\label{TH4}
The codes $VT_0 (n)$ are optimal for all $n$.
\end{conjecture}

This is true for $n \le 8$, as already mentioned, but for larger $n$ it is possible that other, smaller, perfect codes may exist, or even that smaller, optimal but non-perfect codes may exist.

Indeed, consider the code $\{000, 111\}$.
For this code, $\sum_{u \in C} |D_1 (u) | = 1+ 1 =2 < 4$, so this is optimal but not perfect.
For length 4, $\{0000, 0011, 1100, 1111\}$ contains as many codewords as $VT_0 (4)$
(compare \eqn{EqVT1a}), and again is optimal but not perfect.

At length 6 it is possible to
replace two codewords of $VT_0 (6)$ by two other vectors
without affecting its ability to correct single deletions:
110100 and 001011 can be replaced by 111000 and 000111.
The former pair cover eight vectors of length 5, but the latter only
cover four vectors of length 5,
leaving four vectors uncovered.
This suggests the possibility that in some larger code $VT_0 (n)$
it may be possible to replace $k$ vectors by $k+1$ vectors,
which would prove that these codes are not optimal.

In view of these remarks, Conjecture \ref{TH4} does not seem especially compelling!

\section*{Linearity}
As can be seen from \eqn{EqVT1a}, the codes $VT_0 (n)$ are linear for $n \le 4$.
They are never again linear, since, for $n \ge 5$,
$VT_0 (n)$ contains the vectors $1\,0\,0\,0 \ldots 0\,0\,1$ and
$1\,1\,0\,0 \ldots 1\,0\,0$ but not their sum.

In particular, even though $|VT_0 (7) | = 16$, this code is not linear.
One might wonder if it is possible to find a linear code that will do as well, but a computer search has shown that no such code exists.

On the other hand, by adapting a construction of Tenengolts \cite{Ten76},
one can modify the Varshamov-Tenengolts construction to obtain linear codes, with only a small increase in the length of the code.

\begin{definition}\label{de22}
Given $k \ge 1$, let
$$n = k + \left\lceil \sqrt{2k+9/4} + 1/2 \right\rceil ~.$$
The linear single-deletion-correcting code $VT'_0 (n)$ has dimension $k$ and
consists of all vectors $(x_1, \ldots, x_n ) \in \FF_2^n$, where
$x_1, \ldots, x_k$ are information symbols and the $c= n-k$ check symbols
$x_{k+1} , \ldots, x_n$ are chosen so that $\sum_{i=1}^n ix_i \equiv 0$
$(\bmod~n+1)$.
\end{definition}

The construction works because $c$ is just large enough so that
${\binom{c+1}{2}} \ge n+1$, and so
the sums $\sum_{i=k+1}^n ix_i$ cover $n+1$ consecutive values modulo $n+1$.
We omit the details.

The number of check symbols in these codes is of the order of $\sqrt{2n}$,
compared with $O( \log n )$ for the $VT_0 (n)$ codes.
So we end this section with a final question:
What are the optimal linear single-deletion-correcting codes?

\section{Shift register sequences}\label{sectSR}
As mentioned in Section \ref{sectI},
the entry for sequence
\htmladdnormallink{A16}{http://www.research.att.com/cgi-bin/access.cgi/as/njas/sequences/eisA.cgi?Anum=000016}
in \cite{EIS} indicates that these numbers also arise in the enumeration of shift register sequences \cite{Go67}.
We will show here that indeed this is the same sequence.
But whether this is anything more than a coincidence remains an open question.
Of course there are well-known connections between shift-register sequences and conventional error-correcting codes (cf. \cite{MS77}, Chapter 7), so there should be a deeper explanation.

\begin{figure}[htb]
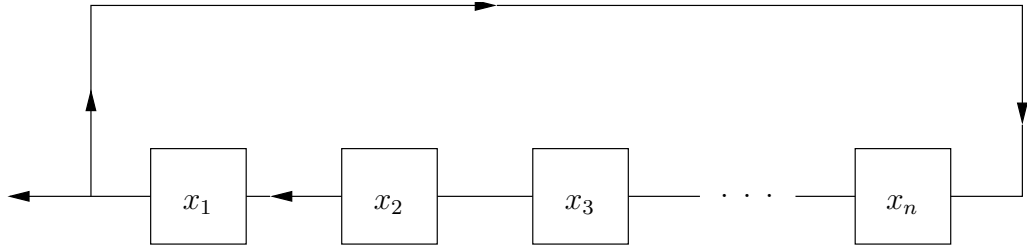

\input difg1.pstex_t
\caption{An $n$-stage pure cycling register.}
\label{FG1}
\end{figure}

The context in which sequence A16 appears in Golomb's book \cite{Go67}
is the enumeration of the (infinite) output sequences from certain types of $n$-stage binary shift registers.
We consider four kinds of shift registers:
the {\em pure cycling register} (or PCR), as illustrated in
Fig.~\ref{FG1},
the {\em complemented cycling register} (or CCR), the {\em pure
summing register} (or PSR) and the {\em complemented summing register} (or CSR).
If the shift register has $n$ cells,
initially containing $x_1, x_2, \ldots, x_n$ $(x_i = 0$ or 1),
then $x_1$ is appended to the output stream, symbols
$x_2, \ldots, x_n$ move to the left, and the symbol
$$
\begin{array}{ll}
\mbox{(PCR)} & x_1 \\
\mbox{(CCR)} & 1+ x_1 \\
\mbox{(PCR)} & x_1 + x_2 + \cdots + x_n \qquad \mbox{or} \\
\mbox{(CSR)} & 1+x_1 + x_2 + \cdots + x_n
\end{array}
$$
is fed back to the right-most cell.

The problem is to determine the numbers of different possible output sequences
from these registers, which we denote by $Z(n)$, $Z^\ast (n)$, $S(n)$ and $S^\ast (n)$,
respectively.
For example $S^\ast (5)=6$, corresponding to the sequences
$$
\begin{array}{l}
\ldots \, 0 \, 0 \, 0 \, 0 \, 0 \, 1 \, 0 \, 0 \, 0 \, 0 \, 0 \, 1 \, \ldots \\
\ldots \, 0 \, 0 \, 0 \, 1 \, 1 \, 1 \, 0 \, 0 \, 0 \, 1 \, 1 \, 1 \, \ldots \\
\ldots \, 0 \, 0 \, 1 \, 0 \, 1 \, 1 \, 0 \, 0 \, 1 \, 0 \, 1 \, 1 \, \ldots \\
\ldots \, 0 \, 1 \, 0 \, 0 \, 1 \, 1 \, 0 \, 1 \, 0 \, 0 \, 1 \, 1 \, \ldots \\
\ldots \, 0 \, 1 \, 0 \, 1 \, 0 \, 1 \, 0 \, 1 \, 0 \, 1 \, 0 \, 1 \, \ldots \\
\ldots \, 0 \, 1 \, 1 \, 1 \, 1 \, 1 \, 0 \, 1 \, 1 \, 1 \, 1 \, 1 \, \ldots ~,
\end{array}
$$
all having period 6 (or a divisor of 6).

Table \ref{TG1}, based on \cite[page 172]{Go67}, shows the first few values of these functions, together with the corresponding sequence numbers
from \cite{EIS}.
\begin{table}[htb]
\caption{Number of ouput sequences from $n$-stage shift registers of types PCR, CCR, PSR, CSR.}

\begin{center}
\begin{tabular}{|c|c|c|c|c|} \hline
~ & PCR & CCR & PSR & CSR \\ \hline
$n$ & $Z(n)$ & $Z^\ast (n)$ & $S(n)$ & $S^\ast (n)$ \\ \hline
1 & 2 & 1 & 2 & 1 \\
2 & 3 & 1 & 2 & 2 \\
3 & 4 & 2 & 4 & 2 \\
4 & 6 & 2 & 4 & 4 \\
5 & 8 & 4 & 8 & 6 \\
6 & 14 & 6 & 10 & 10 \\
7 & 20 & 10 & 20 & 16 \\
8 & 36 & 16 & 30 & 30 \\
9 & 60 & 30 & 56 & 52 \\
10 & 108 & 52 & 94 & 94 \\
$\cdots$ & $\cdots$ & $\cdots$ & $\cdots$ & $\cdots$ \\ \hline
Sequence: &
\htmladdnormallink{A31}{http://www.research.att.com/cgi-bin/access.cgi/as/njas/sequences/eisA.cgi?Anum=000031}
&
\htmladdnormallink{A16}{http://www.research.att.com/cgi-bin/access.cgi/as/njas/sequences/eisA.cgi?Anum=000016}
&
\htmladdnormallink{A13}{http://www.research.att.com/cgi-bin/access.cgi/as/njas/sequences/eisA.cgi?Anum=000013}
&
\htmladdnormallink{A16}{http://www.research.att.com/cgi-bin/access.cgi/as/njas/sequences/eisA.cgi?Anum=000016}
\\ \hline
\end{tabular}
\end{center}
\label{TG1}
\end{table}

Explicit formulas for these functions are given in the next theorem.

\begin{theorem}\label{THSR1}
For $n \ge 1$,
\begin{eqnarray}\label{EqG1}
Z(n) & = & \frac{1}{n} \sum_{d|n} \phi (d) 2^{n/d} ~, \\
\label{EqG2}
Z^\ast (n) & = & S^\ast (n-1) = \frac{1}{2n}
\sum_{d|n \atop d~{\rm odd}} \phi (d) 2^{n/d} \,, \\
\label{EqG3}
S(n) & = & \frac{1}{2(n+1)}
\sum_{d|n+1} \phi (2d) 2^{(n+1)/d} ~.
\end{eqnarray}
\end{theorem}

\begin{remark}\label{rm31}
{\rm Golomb proves \eqn{EqG1} and sketches proofs of the other results.
Actually \eqn{EqG3} is due to Michael Somos (personal communication),
Golomb's version (given in \eqn{EqG4} below)
being slightly more complicated.
The numbers $Z(n)$ (sequence
\htmladdnormallink{A31}{http://www.research.att.com/cgi-bin/access.cgi/as/njas/sequences/eisA.cgi?Anum=000031}) in the
first column are also familiar as the number of binary irreducible polynomials
of degree dividing $n$,
and the number of $n$-bead necklaces formed with beads of two colors, when the necklaces may not be turned over (cf.
\cite[Chap.~4]{Be68},
\cite{GR61},
\cite[Chap.~4]{MS77},
\cite[Problem~7.112]{St99}).
Fredricksen \cite{Fr70} shows that $Z(n) -1$ is the number of 1's in the truth table defining the lexicographically least de~Bruijn cycle.}
\end{remark}

\begin{proof}
Note that sequence
\htmladdnormallink{A16}{http://www.research.att.com/cgi-bin/access.cgi/as/njas/sequences/eisA.cgi?Anum=000016}
appears in two places in the table,
for CCR registers of length $n$ and CSR registers of length $n-1$.
We begin by explaining this, and thus proving that
\begin{equation}\label{EqG3a}
Z^\ast (n) = S^\ast (n-1) ~.
\end{equation}
Suppose for concreteness that $n=4$.
The output sequences from the four types of register are
(omitting plus signs, and writing $1a$ rather than $1+a$, etc.):
$$
\begin{array}{lccccccccccccccc}
{\rm (i)} & a & b & c & d & a & b & c & d & a & \cdots \\
{\rm (ii)} & a & b & c & d & 1a & 1b & 1c & 1d & a & b & c & d & \cdots \\
{\rm (iii)} & a & b & c & d & abcd & a & b & c & d & abcd & a & b & c & d & \cdots \\
{\rm (iv)} & a & b & c & d & 1abcd & a & b & c & d & 1abcd & a & b & c & d & \cdots
\end{array}
$$
In general these sequences have periods $n$, $2n$, $n+1$ and $n+1$, respectively.
If we replace (ii) by the sums of adjacent pairs we get
$$ab \quad bc \quad cd \quad 1ad \quad ab \quad bc \quad cd \quad 1ad \quad
\ldots~,
$$
a CSR(3) sequence.
Conversely, given a CSR(3) sequence, say
$$A \quad B \quad C \quad 1ABC \quad A \quad B \quad C \quad 1ABC \quad \ldots ~,
$$
of period 4, there is a unique CCR(4) sequence of period 8 corresponding to it,
namely
$$0\quad A\quad AB\quad ABC \quad 1 \quad 1A \quad 1AB \quad 1ABC \quad 0 \quad A \quad \ldots \,.
$$
Applying this argument in the general case establishes \eqn{EqG3a}.

In the rest of the proof we make use of Burnside's lemma (cf. \cite{St99}), which states that the number of orbits of a finite permutation group $G$ is equal to the average number of points that are fixed by the elements of $G$.

Let us first prove \eqn{EqG1}.
(This is Golomb's proof \cite[p.~121]{Go67}.)
We take $G$ to be the cyclic group of order $n$ generated by $\pi = (1,2, \ldots, n )$, acting on $\FF_2^n$.
The permutation $\pi^i$ $(1 \le i \le n)$ contains $gcd (n,i)$ cycles,
each of length $n/gcd (n,i)$, and has order $n/gcd (n,i)$.
There are precisely $2^{gcd (n,i)}$ vectors fixed by $\pi^i$,
since each cycle must consist of all 0's or all 1's.
Hence, by Burnside's lemma,
\begin{eqnarray*}
Z(n) & = & \frac{1}{n} \sum_{i=1}^n 2^{gcd (n,i)} \\
& = & \frac{1}{n} \sum_{k|n} \sum_{i=1 \atop gcd (n,i) =k}^n 2^k \\
& = & \frac{1}{n}
\sum_{k|n} 2^k \sum_{gcd \left( \frac{n}{k},i\right) =1} 1 \\
& = & \frac{1}{n} \sum_{k|n} \phi \left( \frac{n}{k}\right) 2^k \\
& = & \frac{1}{n} \sum_{d|n} \phi (d) 2^{n/d} ~.
\end{eqnarray*}

To establish \eqn{EqG2}, we note from (iv) that $S^\ast (n-1)$ is equal to the number of orbits of the same group, but now acting on binary vectors of length $n$ and odd weight.
The number of odd weight vectors fixed by $\pi^i$ is $2^{gcd (n,i) -1}$ if the cycle lengths $n/gcd (n,i)$ are odd, and zero otherwise.
Hence
\begin{eqnarray*}
S^\ast (n-1) & = & \frac{1}{n} \sum_{i=1 \atop n/gcd(n,i)~{\rm odd}}^n
2^{gcd (n,i) -1} \\
& = & \frac{1}{n} \sum_{k|n \atop n/k ~{\rm odd}} \phi
\left( \frac{n}{k} \right) 2^{k-1}
\\
& = & \frac{1}{2n} \sum_{d|n \atop d~{\rm odd}} \phi (d) 2^d ~.
\end{eqnarray*}

Finally, we prove \eqn{EqG3}, by determining $S(n-1)$.
The group is the same, but now (see (iii)) acting on even weight vectors.
If $d= n/gcd (n,i)$ is even there are $2^d$ fixed vectors,
but if $d$ is odd only $2^{d-1}$ fixed vectors.
Hence
\begin{eqnarray}\label{EqG4}
S(n-1) & = & \frac{1}{n} \sum_{d|n \atop d~{\rm odd}} \phi (d)
2^{d-1} + \frac{1}{n} \sum_{d|n \atop d~{\rm even}} \phi (d) 2^d  \\
& = & \frac{1}{2n} \sum_{d|n} \phi (2d) 2^{n/d} ~, \nonumber
\end{eqnarray}
since $\phi (2d) = \phi (d)$ if $d$ odd, $\phi (d) = 2 \phi (d)$ if $d$ even.
\end{proof}

But a mystery still remains: is the fact that the number of codewords in $VT_0 (n)$ equals $Z(n)$ just a numerical coincidence, or is there a one-to-one
correspondence between the codewords and the CCR shift register sequences?
(This is essentially equivalent to a research problem
stated by Stanley in \cite{St86}, Chapter 1, Problem 27(c).)

Furthermore, why is $|VT_1 (n) |$ (sequence
\htmladdnormallink{A48}{http://www.research.att.com/cgi-bin/access.cgi/as/njas/sequences/eisA.cgi?Anum=000048}
in \cite{EIS}), equal to the number of $(n+1)$-bead necklaces with beads of two colors and primitive period $n+1$, when the two colors may be interchanged
but the necklaces may not be turned over (cf. \cite{Fi58}, \cite{GR61})?
This is also the number of irreducible polynomials over $\FF_2$ of degree
$n+1$ in which the coefficient of $x^n$ is 1
\cite{Car52}, \cite{CMRSS}.

\section{Locally transitive tournaments}\label{sectT}
The entry for
\htmladdnormallink{A16}{http://www.research.att.com/cgi-bin/access.cgi/as/njas/sequences/eisA.cgi?Anum=000016}
in \cite{EIS} also indicates that this sequence arose in Brouwer's enumeration \cite{Br80}
of locally transitive tournaments.
A tournament is a directed graph with one directed edge between any two nodes.
It is transitive if there are no directed cycles.

A locally transitive tournament is a tournament such that the subgraphs
on the predecessors of a point and the successors of a point are both transitive.

Brouwer, answering a question raised by P. J. Cameron, determined the number of
locally transitive tournaments on $n$ nodes.
He began by calculating the first few values by computer.
Then he looked up this sequence in 
\cite{HIS}, and found the reference to Golomb's book \cite{Go67}.
With this hint alone, and without having access to the book, he established a
one-to-one correspondence between these tournaments and output sequences from shift registers of CCR type.
From this he obtained the formula
\begin{equation}\label{EqG5}
\sum_{d|n} \,{\rm odd} \left( \frac{n}{d} \right)
\frac{2^{d-1}}{d} \sum_{e \left| \frac{n}{d} \right.} \frac{\mu(e)}{e} ~,
\end{equation}
where odd(i) is 0 or 1 according to whether $i$ is even or odd, and $\mu$ is the
M\"{o}bius function
(\htmladdnormallink{A8683}{http://www.research.att.com/cgi-bin/access.cgi/as/njas/sequences/eisA.cgi?Anum=008683}).
Using the identity
$$\phi (n) = \sum_{d|n} \mu(d) \frac{n}{d}$$
(\cite{Ap76}, p.~26), \eqn{EqG5} immediately reduces to \eqn{EqG2}.

Again we can ask, is there a connection between locally transitive tournaments and the $VT_0 (n)$ codes?

\section{The number of descendants of a vector}\label{sectM}
It was already mentioned in Section \ref{sectVT} that $|D_1 (u) | = r(u)$, the number of runs in $u$.

The next theorem was discovered by E. M. Rains and the author.
Although this must be well-known, we have not found it in the literature.

The {\em derivative} $u' \in \FF_2^{n-1}$ of $u = (u_1, \ldots, u_n) \in \FF_2^n$ is given by
$$u' = (u_1+ u_2, u_2 + u_3, \ldots, u_{n-1} + u_n) ~.
$$
Note that $wt(u') = r(u) -1$.

\begin{theorem}\label{THG1}
\begin{equation}\label{EqGE0}
|D_2 (u) | = {\binom{r(u)+1}{2}} - \delta ~,
\end{equation}
where $\delta = 2 wt(u') - wt (u '' )$ is the {\em deficiency} of $u$.
\end{theorem}

\begin{proof}[Sketch of proof]
First, suppose $u$ is a ``normal'' vector, meaning that all runs have length
$\ge 2$, for example
\begin{equation}\label{EqGE1}
\begin{array}{llcccccccccccccccccccccc}
u & = & 0 & ~ & 0 & ~ & 0 & ~ & 0 & ~ & 1 & ~ & 1 & ~ & 1 & ~ & 0 & ~ & 0 & ~ & 0 \\ 
u' & = & ~ & 0 & ~ & 0 & ~ & 0 & ~ & 1 & ~ & 0 & ~ & 0 & ~ & 1 & ~ & 0 & ~ & 0 \\
u'' & = & ~ & ~ & 0 & ~ & 0 & ~ & 1 & ~ & 1 & ~ & 0 & ~ & 1 & ~ & 1 & ~ & 0
\end{array}
\end{equation}
Then $|D_2 (u) | = \left( \begin{array}{c}
r(u) +1 \\ 2
\end{array}
\right)$ is the number of ways of choosing two things out of $r(u)$ with
repetitions allowed.
If the runs in $u$ have lengths $i,j,k,l, \ldots$, the runs in the shortened vector
have lengths
\begin{equation}\label{EqGE2}
\begin{array}{l}
i-2,~j,~k,~l,~\ldots \\
i,~j-2,~k,~l,~\ldots \\
i,~j,~k-2,~l,~\ldots \\
\quad \cdots~\cdots \\
i-1,~j-1,~k,~l,~\ldots \\
i-1,~j,~k-1,~l,~\ldots \\
\quad \cdots~\cdots
\end{array}
\end{equation}
For a normal vector $wt (u'' ) = 2 wt (u')$ (cf. \eqn{EqGE1}),
$\delta =0$ and \eqn{EqGE0} holds.

Next suppose that all runs in $u$ have length $\ge 2$ except for a single
internal run of length 1, as in
$$
\begin{array}{llccccccccccccccc}
u & = & 0 & ~ & 0 & ~ & 0 & ~ & 0 & ~ & 1 & ~ & 0 & ~ & 0 & ~ & 0 \\
u' & = & ~ & 0 & ~ & 0 & ~ & 0 & ~ & 1 & ~ & 1 & ~ & 0 & ~ & 0 \\
u'' & = & ~ & ~ & 0 & ~ & 0 & ~ & 1 & ~ & 0 & ~ & 1 & ~ & 0
\end{array}
$$
Then $\delta =2$, and indeed $|D_2 (u) |$ is 2 less than it would be for a normal vector, since one of the possibilities in \eqn{EqGE2} vanishes and two others coalesce.

The remaining cases, when there are several runs of length 1, possibly including beginning or ending runs, are left to the reader.
\end{proof}

It is not clear how to generalize Theorem \ref{THG1} to $k$-th order descendants.
Certainly $D_3 (u)$ is not simply a function of the weights of
$u$, $u'$, $u''$ and $u'''$.

\begin{theorem}\label{CHL}
Let
$$\mu_k (n) = \max_{u \in \FF_2^n} |D_k (u) |$$
be the maximal number of $k$-th order descendants of any binary vector
of length $n$. Then
\begin{equation}\label{EqGE5}
\mu_k (n) = \sum_{i=0}^k
{\binom{n-k}{i}} ~,
\end{equation}
for $n \ge k+1$.
Equality is achieved just by the vectors
\begin{equation}\label{EqGE3}
010101 \ldots \quad\mbox{and}\quad 101010 \ldots ~.
\end{equation}
\end{theorem}

According to Calabi and Hartnett \cite{CH69}, \eqn{EqGE5} is proved in an unpublished 1967
report\footnote{I have been unable to locate a copy of this report.}
of Calabi \cite{Cal67}.
The first published proof seems to have been given by
Levenshtein \cite{Lev96}. It was generalized to the nonbinary case by
Hirschberg \cite{Hir99} (see also 
Levenshtein \cite{Lev01} and Hirschberg and Regnier \cite{HR01}).

It is not difficult to show that the vectors \eqn{EqGE3} achieve the bound in \eqn{EqGE5}.

\begin{theorem}\label{THG4}
For the two vectors $010101 \ldots$ and $101010 \ldots$ we have
\begin{equation}\label{EqGE6}
|D_k (u) | = \sum_{i=0}^k {\binom{n-k}{i}} ~.
\end{equation}
\end{theorem}

\begin{proof}
Let $u = 010101 \ldots \in \FF_2^n$, let $M_{n,k}$ be the set of $k$-th order descendants of $u$, and let
$m_{n,k} = | M_{n,k} |$.
Then
\begin{eqnarray}\label{EqGE7}
M_{n,k} & = & 0 | \bar{M}_{n-1, k} \cup M_{n-1, k-1} \nonumber \\
& = & 0 | \bar{M}_{n-1,k} \cup 1 | M_{n-2,k-1} \cup M_{n-2, k-2} ~,
\end{eqnarray}
where the bars denote binary complementation.
However, the last term in \eqn{EqGE7} can be dropped because it is contained in the union of the other two terms.
Since these two terms are disjoint, we have
$$M_{n,k} = M_{n-1,k} + M_{n-2, k-1} ~.$$
This is a disguised version of the recurrence for binomial coefficients, whose solution is given by \eqn{EqGE6}.
\end{proof}

The case $k=2$ of \eqn{EqGE5} is a corollary of Theorem \ref{THG1}:
\begin{corollary}\label{THG2}
For $n \ge 3$,
\begin{equation}\label{EqGE4}
\mu_2 (n) = \sum_{i=0}^2 {\binom{n-2}{i}} =
\frac{1}{2} (n^2 - 3n +4) ~.
\end{equation}
\end{corollary}

\begin{proof}
Let $u$ achieve $\mu_2 (n)$.
The result is easily verified if $r(u)$ is 1 or 2, so
we assume $r(u) \ge 3$.

Suppose $u$ begins with a string of $k \ge 0$ runs of length 1, followed by a run of length $\ge 2$
from position $k+1$.
We will show that the vector $u^\ast$ obtained by complementing $u$ from position $k+2$ onwards
satisfies
$|D_2 (u^\ast ) | \ge |D_2 (u)|$.
By repeating this operation we eventually arrive at one of the vectors \eqn{EqGE3}.

Since $|D_k (u)| = |D_k ( \bar{u} )|$, we may assume that the run following the initial $k$ runs of length 1 in $u$ begins
$1 1 x \ldots$.
In $u^\ast$ this is replaced by $10 \bar{x} \ldots$.
Then we find that $u^\ast$ has $r(u) +1$ runs,
and $wt (u^\ast {''}) = wt (u'' ) -2 + 2x$, from which it follows using
\eqn{EqGE0} that
$|D_2 (u^\ast )| - |D_2 (u) | = r(u) + 2x -3 \ge 0$, as required.
\end{proof}

\section{Related work}\label{sectH}
The history of deletion-correcting codes is closely tied up with studies of codes for correcting other classes of errors such as:
\begin{itemize}
\item
erasures, when bits whose positions are known are deleted
\item
insertions of bits (rather than deletions)
\item
asymmetric errors, when the only errors that occur are that 1's may be changed to 0's
(this is also known as a Z-channel)
\item
unidirectional errors:
0's may be changed to 1's or 1's to 0's, but only one type of error occurs in any particular transmission
\item
bit reversals:
0's may be changed to 1's or vice versa --- this is the subject of classical coding theory
\item
transpositions:
adjacent bits may be swapped
\item
any meaningful combination of the above.
\end{itemize}

Furthermore the alphabet may be changed from $\FF_2$ to $\FF_q$.
This produces an extensive list of families of codes, and of course in each case one can ask
for the largest codes.

In this section we give a brief overview of some other relevant papers.
First, Levenshtein's papers \cite{Lev65}, \cite{Lev65a},
\cite{Lev92}, \cite{Lev01} should be considered essential reading.

Hartnett \cite{Ha74} (see especially Calabi and Hartnett \cite{CH69}) contains some general
investigations of all the above-mentioned codes (both block codes and variable length codes)
from a fairly abstract mathematical point of view.

One of the earliest papers to study deletion-correcting codes is Sellers \cite{Se62}, which combines a special separating string between blocks with a burst-error correcting code inside the blocks.

Ullman \cite{Ull66} uses a construction similar to that of Varshamov and Tenegolts, but his codes
are not as efficient and also use a separating string between blocks.
In \cite{Ull67} he gives bounds on the size of codes for correcting synchronization errors.

Tenengolts \cite{Ten84} generalizes the $VT_a (n)$ codes to larger alphabets.
Nonbinary codes are also discussed in \cite{Bo94},
\cite{Bo95}, \cite{Do85}, \cite{Ma98}.

Other constructions for deletion-correcting and related codes are given by
Calabi and Hartnett \cite{CH69a},
Iizuka, Kasahara and Namekawa \cite{IKN},
Kl{\o}ve \cite{Kl95} and
Tanaka and Kasai \cite{TK76}.

The most recent paper on this subject is by
Schulman and Zuckerman \cite{SZ99}, who present what they describe as ``simple, polynomial-time
encodable and decodable codes which are asymptotically good for channels allowing insertions, deletions and transpositions''.
The number of errors that can be corrected is some constant fraction of the block-length $n$.
The constructions are not explicit.

We conclude this section by mentioning some papers on peripherally related codes.
Codes for correcting asymmetric and unidirectional errors are discussed in
\cite{BR82}, \cite{Et91}, \cite{EO98}, \cite{WVB88} and \cite{WVB89}.
Erasure correcting codes are discussed by Alon and Luby \cite{AL96} and Barg \cite{Ba98}.

\section*{Acknowledgements}
I would like to thank Andries Brouwer, Suhas Diggavi, Vladimir Levenshtein,
Andrew Odlyzko, Eric Rains
and Richard Stanley
for conversations about the subject of this paper;
and David Applegate, David Johnson and Mauricio Resende for their help in establishing that the $VT_0 (n)$ codes are optimal for $n \le 9$ and for their
(so far unsuccessful!) attempts to find better codes.


%

\end{document}